\documentclass[10pt,reqno]{amsart}
\usepackage[T1]{fontenc}
\usepackage{lmodern}
\usepackage[margin=1in]{geometry}
\usepackage{amssymb}
\usepackage[hidelinks]{hyperref}

\newtheorem*{theorem}{Theorem}
\newtheorem{lemma}{Lemma}
\newcommand{\Can}{\mathrm{Can}}
\newcommand{\SCan}{\mathrm{SCan}}
\newcommand{\can}{\mathrm{can}}
\newcommand{\Rel}{\mathrm{Rel}}
\newcommand{\ncl}[1]{\langle\!\langle #1\rangle\!\rangle}

\title[Free Burnside groups of odd exponent]{Free Burnside groups of odd exponent: from 557 to 491}
\author{Elias Judin}
\address[Elias Judin]{University of Cape Town, Cape Town, South Africa}
\email{ejudin@gmail.com}
\author{Stefano Rocca}
\address[Stefano Rocca]{Padova, Italy}
\email{stefano.rocca.99@gmail.com}
\hypersetup{pdfauthor={Elias Judin and Stefano Rocca},pdftitle={Free Burnside groups of odd exponent: from 557 to 491},pdfkeywords={Burnside problem, periodic groups, canonical forms, Lean 4}}

\begin{document}
\begin{abstract}
We prove that the free Burnside group $B(m,n)$ is infinite for every $m\geq 2$ and every odd $n\geq 491$, improving the bound $557$ of Atkarskaya, Rips, and Tent. In their construction of canonical forms, we use the positions of the two ends of a periodic occurrence modulo the length of its primitive root to bound the combined loss under turns of its neighbours. The Fine--Wilf theorem shows that exceeding this bound would contradict the primitivity of the root. We order the turns so that the visited occurrences form an interval and each unvisited occurrence has at most one changed neighbour. In the descent argument, we compare corresponding cyclic cuts before and after turning and account for the common parts of the periodic occurrences. A Lean 4 formalization is given in BurnsideX.
\end{abstract}
\maketitle

Let $F_m$ be a free group of rank $m\geq2$ and put $B(m,n)=F_m/\ncl{w^n:w\in F_m}$. Atkarskaya, Rips, and Tent~\cite[Theorem 2.1]{art} proved infinitude for every odd $n\geq 557$. We prove:

\begin{theorem}
The free Burnside group $B(m,n)$ is infinite for every $m\geq 2$ and every odd $n\geq 491$.
\end{theorem}

We use the notation and the construction by rank of~\cite{art}. Throughout the induction we take $m\geq2$ and odd $n\geq491$. At stage $r\geq1$ we assume IH 1--13 of \cite[Section 3]{art} at all ranks less than $r$, with the estimates at those ranks as modified below. Let $\Rel_r$ be defined as in \cite[Section 4.3]{art}, put $G_r=F_m/\ncl{\bigcup_{j\leq r}\Rel_j}$, and write $\Can_r$ for the set of canonical words of rank $r$. We construct $\can_r(A)$ as the representative in $\Can_r$ of the element of $G_r$ represented by the reduced word $A$. We have $\Rel_0=\{1\}$ and $\can_0$ is free reduction. At ranks less than $r$, existence and uniqueness of the representatives follow from IH 1 and 8. For a fractional power $u$ with specified root $a$, $a^n\in\Rel_r$, its measure is $\Lambda_r(u)=|u|/|a|$ \cite[Definition 3.1]{art}. We call a word $\kappa$-bounded at rank $r$ if every rank-$r$ occurrence has measure $\leq\kappa$. The set of such words in $\Can_{r-1}$ is denoted by $\SCan_{\kappa,r}$ \cite[Definition 4.63]{art}. We retain the constants
\[
 \tau=15,\qquad \mu=n-123,\qquad \varepsilon=2\tau+1=31,\qquad \beta=3\tau+2=47.
\]
For $A=LuR\in\Can_{r-1}$, with $u$ a maximal rank-$r$ occurrence and $a$ its root aligned at its beginning, put $v=a^{-n}u$ (as a reduced word). The turn of $u$ replaces $A$ by $\can_{r-1}(LvR)$. At rank one the shorter side wins since $n$ is odd. At higher ranks, the choice of a side depends on the neighbouring sides. For $A\in\SCan_{\mu,r}$ the procedure of \cite[Section 5.1]{art} selects a side of each member of the stable sequence. We call the resulting word the winner word; its agreement with $\can_r(A)$ is part of the construction. For a cyclic period $B$, let $d'(B)$ be the maximum over its cyclic shifts of the sum of the measures of maximal rank-$r$ occurrences of measure $\geq\beta$ \cite[Lemma 7.25]{art}.

The bound $557$ in~\cite{art} arises from two estimates. In \cite[Section 5.3]{art} the certification thresholds must satisfy
\[
 n/2+77\leq\lambda_2<\lambda_1\leq n-170,\qquad
 \lambda_1-\lambda_2\geq\varepsilon.
\]
The difference allows for the change caused by turning a neighbour. In Case 1 of the proof of \cite[Lemma 7.25]{art}, the estimated decrease of $d'$ is greater than $n-36\tau-16=n-556$. At $n=557$ both restrictions are satisfied. We improve these estimates by distinguishing the ends affected by a turn. The numbers in the headings give the exponent restrictions on the estimates. At $513$ and $505$ we consider when the original side wins. We present the estimates leading to the bound $491$. We are preparing a full account of the proof and working to lower the bound further.

\section{Turns, products and powers}

\noindent\textbf{555. Choice of sides and periodic powers.} In a word of $\SCan_{\mu,r}$, maximal occurrences of measure $\geq5\tau+3=78$ form a stable sequence \cite[Remark 4.56 and Section 5.1]{art}. Following a member through these turns \cite[Convention 4.57]{art}, write $s_i,t_i$ for the measures of its two sides at a choice $i$ of neighbouring sides. \cite[Lemma 4.25(ii)]{art} bounds what survives of the complement from below, and \cite[Corollary 4.36]{art} bounds its maximal prolongation from above. The change from turning one or two neighbours is bounded in \cite[Lemmas 4.44 and 4.59]{art}. Thus
\begin{equation}\label{eq:side-sum}
 n-30<s_i+t_i<n+30,\qquad |t_i-t_j|<2\varepsilon=62.
\end{equation}
Hence $\max_i s_i+\max_j t_j<n+92$, so one of the two sides has measure $<n/2+46$ for every choice of neighbours. To compare opposite sides at a periodic shift, we distinguish the contributions of the two ends. For a fixed side the measure is the sum of the common middle and the two end contributions, in units of $|a|$. Changing the right neighbour preserves the minimizing choice on the left. At the end of a certification sequence only one change of measure $<31$ remains. At a periodic shift on the boundary, the preceding letter and the common segment from the cut $47$ periods into the occurrence to the next occurrence determine the change at the right end. By \eqref{eq:side-sum} and the one-neighbour estimate, the sum of the measures of the two opposite sides is less than $n+30+3\varepsilon$. If both measures exceeded $n/2+61\frac12$, their sum would exceed this bound. As in the proof of \cite[Lemma 7.10]{art}, agreement along the shifts would force a forbidden overlap $BBb$ in the threshold sequence (see~\cite{fife}), with $b$ the first letter of $B$. At $555$ the thresholds are $354$ and $385$.

In \cite[Proposition 7.19]{art}, $M_3$ in the canonical product must be $\tau$-free modulo rank $r$ \cite[Definition 3.6]{art}. Suppose $C^\tau$ occurs there, with $C$ primitive, cyclically canonical at rank $r-1$, and $C^n\notin\bigcup_{j\leq r}\Rel_j$. We consider first the case where the prepared middle survives and $C^\tau$ has no interior maximal occurrence of measure $\geq78$. If the surviving part of the complement has measure $<31$, the two pieces cut from it each have measure $<15$, giving raw measure $<61$. Subtracting the old middle and the opposite input leaves more than $n-61-46-(n/2+46)=n/2-153>78$ on each input. By \cite[Corollary 4.18]{art}, each such subword lies in a unique maximal occurrence of the indicated rank and root, of measure greater than $78$. Otherwise the surviving part has measure at least $31$. We count the shifts of the rank-$r$ $15$th power in the subword specified in \cite[Remark 7.8]{art}. As in \cite[Lemma 7.9]{art}, at least $\tau-4=11$ remain between the bounding occurrences. At most two cross a cut, at most three lie in the middle, and at most four lie on each side. Hence both sides contain a shift. The surviving part lies in $C^\tau$, so a cyclic shift of $C$ contains a rank-$r$ $30$th power \cite[Remark 7.2]{art}. Counting shifts of this longer subword gives at most $3+2+3+2=10$ (left, middle, right, and crossing), fewer than eleven.

In the descent argument of \cite[Lemma 7.25]{art}, we show that the decrease per period is at least $n-518$. A turn removes an occurrence of measure $>\mu-\tau=n-138$. At most two counted occurrences cross each cut. After matching the unchanged parts, each contributes less than $47$ periods of unmatched length there, and at most one contributes $16$ or more. With a $\tau$th power surviving on each side of the turn, the two free ends and the fractional power between them contribute at most $30+170+30=230$. Including both cuts gives $<230+2(47+16)=356$. We turn corresponding occurrences in $N$ periods. Dividing the decrease by $N$ and letting $N\to\infty$ gives a decrease of at least $n-494$ per period, since the end corrections are bounded independently of $N$. If the retained arc contains no such $\tau$th power, the new measure is at most $380$, giving a common decrease $n-518$.

\medskip
\noindent\textbf{553. Uniqueness of the seam turn.} In \cite[Lemma 4.70]{art}, a product is prepared by successive seam turns, with margin $46$. An intermediate word is $A'QC'$ with $\kappa$-bounded outside words and a $(3\tau+1)$-free middle of rank $r$ (the excluded powers have $a^n\in\Rel_r$). We show that if a maximal occurrence has measure at least $\kappa+47$, no other occurrence permits a turn of measure $\geq\kappa+15$. The former occurrence extends a whole period of length $p$ beyond each end of $Q$, since the middle contributes less than $46$ periods and either outside piece at most $\kappa$. Distinct maximal occurrences cannot contain one another. Let a competing occurrence of period $q$ start before the occurrence of period $p$. It must end before it. If it extends $t>0$ letters past the left cut, their overlap contains at least $p+t$ letters. By \cite[Corollary 4.19]{art}, $p+t<16\min(p,q)$, and hence $t<15q$. The other case follows by symmetry. Thus the occurrence of period $p$ is the only occurrence permitting a turn of measure $\geq\kappa+15$ while it is present. For this comparison the margin $47$ replaces the general margin $62$.

\medskip
\noindent\textbf{551 and 549. Periodic factors across cuts.} We show that a periodic factor through $15$-, $166$-, and $15$-free pieces has measure $<197$. It contains at most $14+165+14$ aligned whole periods inside the pieces and at most two crossing the joins. The two incomplete ends bring its measure to $<197$, within the bound $n/2-77$ of \cite[Lemma 5.12]{art}. In the product with one more letter, we prepare the prefix with bound $n/2+47$. The remaining changes give $n/2+47+46+15\leq\mu$. For an aligned periodic prefix of period length $d$, the $15$-free correction can remove at most $15d-1$ letters and the appended letter at most one. Hence the required subword survives if the prefix exceeds its length by more than $15d$ letters. The remaining cyclic bound is $2\cdot198+15+1=412\leq n-138$, first valid at odd $n=551$.

Let the cyclic period be $XY$, with $a$-periodic factors of measure $<15$ in $X$ and $\leq198$ in $Y$. We show that no repetition of $XY$ contains $a^{411}$ if $a^n\in\Rel_r$ and $(XY)^n\neq1$ in $G_r$. In $XYX$ or $YXY$ the measure of such factors is $<411$. An $a^{411}$ in a repetition of $XY$ must therefore contain both pieces in full, whence $|XY|<213|a|$. Its length exceeds $|XY|+|a|$, so by the Fine--Wilf theorem~\cite[Lemma 4.9]{art}, \cite{finewilf}, $XY$ is conjugate to a power of $a$. Since $a^n\in\Rel_r$, this contradicts $(XY)^n\neq1$ in $G_r$. The inequality $411\leq n-138$ gives $549$.

\medskip
\noindent\textbf{525. Powers meeting correction words.} In the situation of \cite[Lemma 7.24]{art}, write a sufficiently high canonical power as $Dva^NwE$, with $a$ primitive and $a^n\notin\Rel_r$. Here $va^Nw$ is $a$-periodic and $166$-free of rank $r$, and $D,E$ are $15$-free of that rank. We show that for sufficiently large $N$ a rank-$r$ periodic factor of measure $>4$ cannot meet both $D$ and $E$. Let $N$ be such that $T=N|a|\geq\max(2|a|,|Dv|+|wE|)$. A $b$-periodic factor with $b^n\in\Rel_r$ meeting both corrections contains $a^N$ and has length at most $2T$. If it is longer than $4|b|$, then $T>2|b|$ and $T>|a|+|b|$. By \cite[Lemma 4.9]{art}, $a$ would be a cyclic shift of $b$, forcing $a^n\in\Rel_r$ by \cite[IH 3]{art}. Hence a factor of measure $>4$ meets at most one correction, giving the bound $181=166+15$. In the cyclic argument, a fractional power with measure in $(n-139,n+139]$ has a raw complement of at most $139$ periods, giving $169$ after normalization. For a $169$-bounded word with a $45$th power of rank $r$ in its cyclic core, the refinement of \cite[Lemma 7.23]{art} at $549$ gives $2\cdot169+15=353$ for every power of the resulting period. The remaining restriction comes from certification. The thresholds $n/2+61\frac12$ and $n/2+92\frac12$ differ by $31$, and the upper one is at most $n-170$ precisely for $n\geq525$.

\section{Neighbouring turns and periodic shifts}

\noindent\textbf{523. The two ends of a periodic occurrence.} Let $A\in\SCan_{\mu+1,r}$ and enumerate its maximal rank-$r$ occurrences of measure at least $78$ as $u_1,\ldots,u_\ell$ from left to right. For $T\subseteq\{1,\ldots,\ell\}$, let $A_T$ be the word obtained by choosing the complementary sides at the indices in $T$. We extend the construction of the family from the source bound $\mu$ to $\mu+1$. In particular, $A_\varnothing=A$, each $A_T\in\Can_{r-1}$, and changing one index is a turn of the corresponding maximal occurrence. Both sides at an index have the same period length and their measures depend only on the choices at the adjacent indices. The sum and neighbour estimates in \eqref{eq:side-sum} hold, and the occurrences corresponding to the indices have measure at least $16$. These properties do not require $A_T\in\SCan_{\mu,r}$.

\begin{lemma}\label{lem:two-end}
Fix an index of this family and one of its two sides. Let $s_i,s_j$ be its measures at two choices of neighbouring sides. If its length at choice $i$ is $qd+\rho$, where $d$ is its period length and $0\leq\rho<d$, then its length at choice $j$ is greater than $(q-61)d$. In particular,
\begin{equation}\label{eq:two-end-loss}
 s_i\geq k+61\implies s_j>k\qquad(k\in\mathbb N).
\end{equation}
\end{lemma}
\begin{proof}
Put $L=qd+\rho$, and let $L_l,L_r$ be the lengths obtained by making the left and right changes separately, both measured from choice $i$. If a neighbour is unchanged the corresponding length is $L$. For a fixed side the length is a sum of contributions from the middle and the two ends. Thus, with $L'$ the length after both changes, $L+L'=L_l+L_r$. Each of $L-L_l,L-L_r$ is less than $31d$. If either is at most $30d$, then $L-L'<61d$.

Suppose both losses exceed $30d$. Write $15d+x,15d+y$ for the lengths of the two overlaps at choice $i$. The respective losses are less than $30d+x,30d+y$. By same-rank nonnesting a neighbour of period at least $15d$ gives a loss of less than $30d$ letters. Thus the two periods are $15d-u,15d-v$, where $u,v>0$. By \cite[Lemma 4.9]{art}, $x+u<d$ and $y+v<d$.

Reducing the neighbour periods modulo $d$ gives intervals of the $d$-periodic repetition of the root,
\[
 [d-u,d+x+u),\qquad [d+\rho-y-v,d+\rho+v),
\]
with periods $u$ and $v$, respectively. If $x+y\geq d+\rho$, their intersection is $[d-u,d+\rho+v)$ of length $\rho+u+v$. By Fine--Wilf it has period $g=\gcd(u,v)$. It contains at least $u$ and at least $v$ letters. Since $g$ divides both $u$ and $v$, the periodicity of the intersection extends to both intervals: a whole period of either interval is contained in the intersection. Their union has length $x+y+u+v-\rho\geq d+g$. Equality of letters at distance $g$ therefore holds for $d$ consecutive starting positions and hence throughout the $d$-periodic repetition. Since $0<g<d$, the repetition has period $\gcd(d,g)<d$, contrary to primitivity. It follows that $x+y<d+\rho$ and hence $L-L'<61d+\rho$. Thus $L'>(q-61)d$.
\end{proof}

If the opposite side is bounded by $n/2+92$, \eqref{eq:side-sum} gives initial measure $>n/2-122=139\frac12$ at $523$. By Lemma~\ref{lem:two-end}, more than $78$ periods survive, as required in \cite[Section 5.1]{art}.

\medskip
\noindent\textbf{519. Periodic shifts between stopping cuts.} We exclude the case in \cite[Proposition 7.19]{art} where $C^{15}$ contains a maximal rank-$r$ occurrence of measure at least $78$ strictly in its interior. By \cite[Remark 7.2]{art}, there are at least $15-2=13$ interior shifts of this occurrence. We consider the eleven shifts other than the first and last of these thirteen. As in \cite[Lemma 7.10]{art}, let $k$ be the advance in the occurrence index over one period. Between the stopping cuts at most five occurrences are common to the two sequences, at most one on each outside piece and at most three meeting the middle. \par\noindent (a) Suppose a neighbour choice takes one of these eleven occurrences below $78$. By periodicity the same choice applies to all eleven, placing them between the cuts. At most two common occurrences are tested shifts already turned. Each missing shift is adjacent to a common occurrence already turned. For $k\geq2$, each common occurrence accounts for at most two shifts, giving at most $10$. For $k=1$, each accounts for at most one, except that a tested occurrence already turned can account for three. This gives at most $5+2\cdot2=9$. \par\noindent (b) Suppose the eleven occurrences persist. A missing outermost shift would leave its outside scan empty and only $5+4$ places for the eleven. Thus all thirteen persist, split $4+5+4$. Four shifts on each outside scan make the extreme shifts their starting occurrences. Comparing translated positions in the two scans forces the initial side choice to repeat one period later, giving the forbidden overlap in \cite[Lemma 7.10]{art}. At $519$ the cyclic decrease $n-518$ is one.

\medskip
\noindent\textbf{517. Certification and the common endpoint of two cuts.} We raise the upper bound on the certification thresholds to $n-169$. If one choice gives $s_j<78$, \eqref{eq:two-end-loss} gives $s_i<139$ for every choice, and \eqref{eq:side-sum} gives $t_i>n-169$. The thresholds $n/2+58\frac12$, $n/2+89\frac12$ fit at $517$. In the cyclic comparison with a nonempty common arc, both cuts retain its right endpoint. Counting its correction once lowers the bound for the surviving fractional power by $\tau$, from $168$ to $153$. The total measure falls from $380$ to $365$, giving decrease $n-503$.

\medskip
\noindent\textbf{515, 513, and 505. A criterion for the original side to win.} In \cite[Lemma 7.22]{art} the cyclic period $b=EW'DW''$ consists of $15$-, $45$-, $15$-, and $45$-free pieces. We show that its repetitions are $165$-free of rank $r$ unless $b=a^s$ for some $a^n\in\Rel_r$ and integer $s\geq0$. Suppose $c^{165}$ occurs in a repetition of $b$, where $c^n\in\Rel_r$. Starting at the first letter of $c^{165}$, we read one cycle and the remainder of the piece split by the starting cut. This interval excludes $120+45=165$ whole copies, since only the split piece is read twice. It contains a full cycle starting at a piece boundary, which excludes $120$ copies. The occurrence of $c^{165}$ extends beyond the interval and contains that cycle. Thus $|b|<121|c|$, since $121|c|$ consecutive letters of a $c$-power contain $c^{120}$. By \cite[Lemma 4.9]{art}, the periodic repetitions have a common primitive root up to cyclic shift; hence $b=a^s$ with $a^n\in\Rel_r$. Otherwise the middle is $165$-free, and the one-end bound becomes $180$, within $n/2-77$ at $515$.

For an original side of measure $\leq k$, the contrapositive of \eqref{eq:two-end-loss} gives $s_i<k+61$ for all $i$. Every complement has measure $>n-30-(k+61)$. Thus the original side wins whenever $2(k+61)\leq n-30$, or $2k+153\leq n$ for odd $n$. Taking $k=180$ gives $513$. At $505$, only the first and last members of the sequence can meet the correction ends. They have at most one close neighbour, so their original sides stay below $180+31=211<(505-30)/2$. Interior members have measure $<166$ in the $165$-free middle, hence satisfy the same criterion with $k=176=(505-153)/2$.

\section{Canonical forms and cyclic complexity}

\noindent\textbf{495. Measures of occurrences in canonical forms.} For odd $n\geq557$ and $A\in\SCan_{\mu,r}$, \cite[Lemma 5.19]{art} gives $\can_r(A)\in\SCan_{n/2+46,r}$. Let $s$ denote the winner side and $t$ the other side. The minimum measure of $s$ over certified neighbour choices is at most the corresponding minimum for $t$. The neighbouring winners give a choice $i$ admissible for both \cite[Lemma 5.17]{art}. Hence
\[
 s_i<\min_j s_j+62\leq t_i+62,\qquad s_i+t_i<n+30.
\]
Combining the inequalities gives $s_i<n/2+46$. Every maximal occurrence of measure at least $140$ in the resulting word comes from the stable sequence (see Lemma~\ref{lem:recognition}(a) below), so the bound holds for the word. For the comparison of semicanonical forms in \cite[Section 6]{art}, we use the successive bounds $n/2+15$, $n/2+46$, and $\kappa=n/2+77$, one neighbour change apart. Preparation must stay inside $\mu$ and the complement of a long turn must be recognized. These conditions are
\begin{equation}\label{eq:recognition}
 \kappa+47\leq\mu,\qquad n-30-\kappa\geq140.
\end{equation}
Each is equivalent to $n\geq494$. In the cyclic case with no surviving $\tau$th power on the retained arc, adjacent free pieces can be estimated together. The resulting measure is at most $350$ for an empty common arc and at most $46+30+(153+32)+30+46+15=352$ otherwise. The larger cost is $356$ in the case with a $\tau$th power on both sides. Thus $d'$ decreases by at least $n-138-356=n-494$, giving $495$.

\medskip
\noindent\textbf{493. Recognition of occurrences and seam turns.} The following assertion refines \cite[Corollary 4.74]{art}.

\begin{lemma}\label{lem:recognition}
Let $A\in\SCan_{\mu+1,r}$ and let $A_T$ be the family constructed at $523$.
\begin{enumerate}
\item[(a)] Every maximal rank-$r$ occurrence in $A_T$ of measure at least $139$ corresponds to an index of the family, with its aligned root.
\item[(b)] If $A=X_0\mapsto X_1\mapsto\cdots\mapsto X_l$ is a chain of rank-$r$ turns of maximal occurrences of measure at least $139$, then each $X_j$ is some $A_{T_j}$. No bound of $\mu$ on the $X_j$ is required.
\end{enumerate}
\end{lemma}
\begin{proof}
Suppose an occurrence $u$ in (a) does not correspond to an index. Its starting position is strictly between two successive labelled starting positions or in an exterior interval. Undoing the turns outside the indices bordering this interval preserves $u$ and its position relative to the labelled occurrences. Each boundary pullback preserves the period length $d$ and loses less than $31d$ letters. If there are two, undo the left turn first. The common suffix has length greater than $|u|-31d\geq108d$. The right overlap has length less than $16d$ and thus is a suffix of the original $u$. Both overlaps have the phase determined by the ends of $u$. If either loss is at most $30d$ their sum is less than $61d$. Otherwise the overlap estimates and the interval argument of Lemma~\ref{lem:two-end} give $x+y<d+\rho$, where $|u|=qd+\rho$. The losses are at most $30d+x$ and $30d+y$, so their sum is less than $61d+\rho$. This application uses no length identity for an occurrence outside the family.

Since $q\geq139$, the maximal occurrence obtained in $A$ has measure greater than $q-61\geq78$; with at most one boundary pullback the loss is less than $31d$, giving the same conclusion. It must belong to the enumeration in $A$, contrary to the preservation of its strict position relative to the family. The aligned root is unique by the overlap bound. This proves (a).

For (b), $T_0=\varnothing$. If $X_j=A_{T_j}$, (a) identifies the next turned occurrence with an index $i$. The defining turn property gives $X_{j+1}=A_{T_j\mathbin{\triangle}\{i\}}$, where $\triangle$ denotes symmetric difference. Induction proves (b).
\end{proof}

At the seam, half a period $p$ beyond each end of $Q$ suffices. If $p\geq2q$, that half uses at least $q$ of the fewer than $16q$ letters of overlap. If $p<2q$, Fine--Wilf bounds the overlap of distinct maximal occurrences by $p+q<3q$. In both cases penetration is $<15q$ as required at $553$. The new inequalities are $\kappa+46\frac12\leq\mu$ and $n-30-\kappa\geq139$. Both hold at $493$, where the lower certification threshold is $n/2+46\frac12=n-200$. In the cyclic estimate the left cut and its adjacent correction contribute at most $78$ together. This replaces $63+30$ by $78$ and reduces $356$ to $78+200+63=341$.

\medskip
\noindent\textbf{491. Winner sides and the order of turns.} Let $n=2h+1$ and, as at $495$, let $s$ denote the winner side and $t$ the other side. Let $i$ be the choice of neighbouring winner sides, admissible for both. We show that $s_i<h+46$. If $s_i\geq h+46$, then \eqref{eq:two-end-loss} with $k=h-15$ would give $\min_j s_j>h-15$. The winner rule gives $t_i\geq\min_j t_j\geq\min_j s_j$, contradicting $s_i+t_i<n+30=2h+31$. Thus, for the choice of neighbouring winners,
\begin{equation}\label{eq:integral-winner}
 s_i<h+46.
\end{equation}
By Lemma~\ref{lem:recognition}(a), the winner word is $(h+46)$-bounded. At $491$ this is $291$, half a period below $n/2+46$.

\textit{Certification} \cite[Definition 5.3 and Lemma 7.10]{art}. We show that the periodic comparison holds with the thresholds $\lambda_2=h+46$ and $\lambda_1=h+77$ in the overlap-free sequence of thresholds $\lambda_{m(j)}$ from~\cite{art}, starting with the high value. At $491$ these are $291=n-200$ and $322=n-169$. Two low values sum to $n+91$, one short of $n+30+2\varepsilon$ at a periodic shift where the left neighbours have yet to be matched. Let $k$ count occurrences per period. \par\noindent (a) If $k\geq3$, the comparison starts one occurrence later, where both neighbours match.\par\noindent (b) If $k=1$, two low values force the next one to be high, and the following two fully matched occurrences have the same side. Comparison with the high value requires $n+30+\varepsilon$.\par\noindent (c) If $k=2$, disagreement would force three low values two indices apart. The intervening five values would contain $000$ or $01010$ (low $=0$), both forbidden overlaps. \par In the other boundary comparisons at least one threshold is $\lambda_1$, so their sum is at least $\lambda_1+\lambda_2=n+122>n+30+2\varepsilon$.

\textit{Comparing semicanonical forms} \cite[Section 6]{art}. At $491$, $\kappa+46\frac12=369=\mu+1$. Let $A\in\SCan_{\mu+1,r}$ be the common source of two chains of turns of measure at least $139$, ending at $B$ and $C$, where $C$ is the chosen $(h+46)$-bounded endpoint. By Lemma~\ref{lem:recognition}(b), $B=A_S$ and $C=A_T$ for subsets $S,T$ of the same family. Starting at $A_S$, change the indices of $S\mathbin{\triangle}T$ in increasing order. At each changed index the new side and its left neighbour agree with their choices in $A_T$. The measure of the new side is less than $h+46+31$ by the one-neighbour estimate. By \eqref{eq:side-sum}, the old side has measure greater than
\[
 n-30-(h+46+31)=h-106\geq139.
\]
Thus $B$ is joined to $C$ by turns of measure at least $139$. The margin in the last inequality, in $\kappa+46\frac12\leq n-122$, and in $h+77\leq n-169$ is $h-245\geq0$. The auxiliary source bound is not the domain of invariance of the winner rule.

\textit{Agreement with the canonical form} \cite[Lemma 6.22 and Proposition 5.20]{art}. The prepared-product bound $n/2+92$ exceeds $\mu-31$ by half a period at $491$. At most one original occurrence in the prepared product has measure greater than $n-169$. Indeed, if the word is $(h+77)$-bounded there is none. Otherwise its $16$-free middle lies between two $15$-free corrections and two $(h+46)$-bounded outside words. Since $h+46+15+16\leq n-169$, each occurrence of measure greater than $n-169$ contains sixteen of its periods starting at the first cut of the middle. By \cite[Corollary 4.19]{art}, two such occurrences coincide.

\begin{lemma}\label{lem:ordered-turns}
Let $A=A_\varnothing\in\SCan_{\mu,r}$ have at most one member of its family of measure greater than $n-169$. Let $W=A_T\in\SCan_{h+46,r}$ be its chosen winner word. There is a chain of turns from $A$ to $W$ such that each turn is of a maximal rank-$r$ occurrence of measure at least $16$, and all words of the chain belong to $\SCan_{\mu,r}$.
\end{lemma}
\begin{proof}
If $T=\varnothing$, there is nothing to prove. Choose the index $e$ of the member of measure greater than $n-169$, if there is one, and any index if there is none. Visit $e$, then the indices to its left in decreasing order, and finally those to its right in increasing order. A visit changes the side if and only if its index belongs to $T$. Let $M_j\subseteq T$ be the set of indices changed after $j$ visits. The word is $A_{M_j}$ by the defining turn property. We have $M_0=\varnothing$, $M_\ell=T$, and the visited indices form an interval containing $e$ for $j>0$.

We prove by induction on $j$ that $A_{M_j}\in\SCan_{\mu,r}$. The assertion holds for $j=0$. For $j>0$, consider a maximal occurrence in $A_{M_j}$. There are three cases.

\noindent (a) It corresponds to an index whose side agrees with its side in $W$. This includes every index not in $T$. Its measure differs from its measure in $W$ by less than $62$, and hence is less than $h+46+62\leq\mu$.

\noindent (b) It corresponds to an index in $T\setminus M_j$. This index has not been visited and is different from $e$. At most one neighbour has changed, since the visited indices form an interval. Its measure in $A$ is at most $n-169$, so its measure in $A_{M_j}$ is less than $n-169+31=n-138<\mu$.

\noindent (c) It corresponds to no index. By Lemma~\ref{lem:recognition}(a), its measure is less than $139<\mu$.

Each word of the family is in $\Can_{r-1}$. The three cases therefore give $A_{M_j}\in\SCan_{\mu,r}$ without assuming that bound when applying the recognition lemma. At a visit that changes a side, the family properties ensure that the measure of the turned occurrence is at least $16$. Other visits leave the word unchanged.
\end{proof}

The prepared product has bound $n/2+92\leq\mu$, and its winner satisfies \eqref{eq:integral-winner} and Lemma~\ref{lem:recognition}(a). Thus Lemma~\ref{lem:ordered-turns} applies. The refined winner rule is invariant under turns of measure at least $16$ between words of $\SCan_{\mu,r}$, as in \cite[Proposition 5.20]{art}. Along the chain the winner remains $W$. The rank-$r$ normalization identifies $W$ with the unique canonical representative of its $G_r$-class. Consequently $\can_r(A)=W$, as required for the canonical product. The normalization result precedes this use of canonical agreement in the construction of the triangles.

\textit{Descent and infinitude.} In the descent branch of \cite[Lemma 7.25]{art}, the new period $C$ is conjugate to $B$ in $G_r$. We show that $d'(C)+n-490\leq d'(B)$. For each cyclic cut of $C$ we compare a corresponding cut of $B$. Deleting $n$ whole periods saves at least $n-440$ at these cuts. For the short turns, comparison of coefficients in the number of periods saves at least $n-138-341=n-479$ or $n-138-352=n-490$. Taking the maximum over cuts of $C$ gives
\begin{equation}\label{eq:cyclic-descent}
 d'(C)+n-490\leq d'(B).
\end{equation}
By \eqref{eq:cyclic-descent}, for odd $n\geq491$ the integer part of $d'$ decreases by at least one. The alternative gives a $45$-free cyclic core or a period whose powers are $353$-bounded, with $353\leq n-138$. The inequality becomes an equality at $491$.

The rank-zero and rank-one constructions give the initial ranks. At a positive successor the relator definition is unchanged; closure under inverses and cyclic shifts, disjointness of ranks and the small-cancellation estimates give IH 3 and 5. Certification and semicanonical comparison enter the construction of the canonical set and the unique representative of a $G_r$-class (IH 1, 2, 4, 6--9). For IH 10--12 one also needs the multiplication and periodic gluing statements, including compatibility of the canonical triangles with those of the preceding ranks. The agreement with the winner word and the product and prefix estimates above enter these statements. Given these laws at rank $r$, the descent alternatives, the description of powers and gluing with the conjugating words give IH 13, as in \cite[Proposition 7.26]{art}. The successor construction assembles the results and extends the choices of canonical triangles without changing the preceding ranks. The induction requires these geometric and normalization results as well as the descent estimate.

As in \cite[Proposition 8.6 and Corollary 8.7]{art}, persistence gives $\ncl{\bigcup_r\Rel_r}=\ncl{w^n:w\in F_m}$. If two words are equal in this quotient, their difference is a finite product of conjugates of relators and their inverses. They are therefore equal in some $G_r$, and have equal rank-$r$ canonical forms \cite[IH 8 and Theorem 8.1]{art}. Using two distinct free generators as the letters, the prefixes of an infinite overlap-free binary word~\cite{fife} are reduced, cube-free and fixed by every $\can_r$ \cite[Remark 8.8]{art}. Distinct prefixes therefore give distinct elements of $B(m,n)$.

\section*{Acknowledgments}
We developed BurnsideART~\cite{burnsideart}, a Lean formalization of~\cite{art}, with Codex 6.1 Sol. We developed BurnsideX~\cite{burnsidex}, the Lean formalization accompanying this paper, with Codex 6.2 Sol. ChatGPT assisted with manuscript preparation.

\end{document}